\documentclass[a4paper,12pt]{amsart}
\usepackage[T1]{fontenc} 

\usepackage{graphics}
\usepackage{amsmath}
\usepackage{mathrsfs}
\usepackage{stmaryrd}
\usepackage{eucal}
\title{Asymptotic freeness of Jucys-Murphy element and a certain projection}
\author{Lech Jankowski}

\newtheorem{tw}{Theorem}

\newtheorem{lem}[]{Lemma}
\newtheorem{fa}[]{Fact}
\newtheorem{wn}{Wniosek}

\theoremstyle{definition}
\newtheorem{defi}{Definition}
\newtheorem{pr}[]{Example}

\theoremstyle{remark}
\newtheorem{uw}[tw]{Remark}

\DeclareMathOperator{\tr}{tr}
\DeclareMathOperator{\Tr}{Tr}

\providecommand{\norm}[1]{\left\lVert#1\right\rVert}
\providecommand{\abs}[1]{\left\lvert#1\right\rvert}

\begin{document}

\begin{abstract}
We explain the appearance of the free compression of a transition measure in the problem of the restriction of the representation of the symmetric group
to a subgroup by showing the responsible free projection.
\end{abstract}
\maketitle

\section{Introduction}

It was shown in [Bia98] that the restriction of a representation of the symmetric group $S_n$ to the subgroup $S_k$ (the inclusion is defined by
declaring numbers bigger than $k$ as fixed points) can be descried by free compression of measure. The measure is the Kerov transition measure
of the representation. It is not fully understood as there are no free random variables known to be responsible for this phenomenon.
In this paper we find a projection which is free from Jucys-Murphy element known as a random variable which distribution in a certain 
noncommutative probability space is equal to the transition measure of an arbitrary representation of $S_n$. This gives a conceptual explanation of the
phenomenon discovered by Biane.

For esthetical reason our starting point is a non-commutative probability space $(\mathbb{C}[S_{n+1}],\operatorname{tr}\rho(\bullet\downarrow_{S_n}^{S_{n+1}}))$
and its element $X=(1,n+1)+(2,n+1)+\cdots+(n,n+1)$ called \emph{Jucys-Murphy element}.

In order to prove our result we need to extend this probability space using the following idea from [Bia98, proof of Prop.~3.3]:
`\emph{We identify $S_{n+1}$ with $S_n\times\{e,(1,n+1),(2,n+1),\dots,(n,n+1)\}$ by the map $(\sigma,\tau)\rightarrow\sigma\tau$.
In this way we can represent an operator on $\mathbb{C}[S_{n+1}]$ by an $(n+1)\times(n+1)$ matrix of operators on $\mathbb{C}[S_{n}]$.}'
In the rest of this article we will thus work in a space $(\mathbb{C}[S_{n}]\otimes \operatorname{End}(\mathbb{C}^{n+1})), \operatorname{tr}\rho(\bullet)\otimes \operatorname{tr}(\bullet) )$.

It was shown by Biane [Bia98, proof of Prop.~3.3] that the action of $X$ by the left regular representation of $\mathbb{C}[S_{n+1}]$ is 
represented by a matrix 
\begin{equation} \label{eq:macierzjucysia}
X=\left[
\begin{matrix}  
 0& 1 & 1 & 1 & \dots & 1 & 1\\
 1& 0 & (1,2) & (1,3) & \dots & (1,q-1)& (1,q)  \\
 1& (1,2) & 0 & (2,3) & \dots  & (2,q-1) & (2,q) \\
 1& (1,3) & (2,3) & 0  & \dots  & (3,q-1) & (3,q)\\
 \vdots & \vdots & \vdots   & \vdots & \ddots & \vdots & \vdots\\
 1& (1,q) & (2,q) & (3,q) & \dots  & (q-1,q) & 0 \\
 \end{matrix} \right],
 \end{equation}
 where the entries are acting by the right regular representation of $\mathbb{C}[S_{n}]$ and that the distribution of $X$ is the transition measure of
 the representation $\rho$.

\section{The result}

We are interested in restricting representations from $S_n$ to $S_k$, thus we are looking for a projection which can compress $J_{n+1}$ to 
$J_k$. Clearly, $P$ given by a matrix 
\begin{equation} \label{eq:macierzjucysia}
P=\left[
\begin{matrix}   1 & 0 & 0 & \dots & 0 & 0 \\
 0 & 1 & 0 & \dots  & 0 & 0 \\
0 & 0 & 1 & \dots  & 0 & 0 \\
 \vdots & \vdots & \vdots & \ddots  & \vdots & \vdots \\
 0 & 0 & 0 & \dots  & 0 & 0 \\
 0 & 0 & 0 & \dots & 0 & 0
 \end{matrix} \right],
 \end{equation}
 where $1$ occurs $k+1$ times
 has the desired property. Let us describe this projection in the language of
the group algebra $\mathbb{C}[S_{n+1}]$. 
In order to do that we will evaluate $P$ on the elements of the basis. For $\sigma\in S_n$ we have 
$$
P(\sigma)=\sigma,
$$
$$
P(\sigma(j,q+1))=\sigma(j,q+1),\,\text{for all } j\in\{1,2,\dots,k\},
$$
$$
P(\sigma(j,q+1))=0,\,\text{ for all } j\in\{k+1,k+2,\dots,q\}.
$$
One can conclude that $P(\tau) = \begin{cases} \tau &\mbox{if } \tau^{-1}(q+1)\in\{1,2,\dots,k,q+1\} \\
0 & \mbox{if } \tau^{-1}(q+1)\in\{k+1,k+2,\dots,q \}.  \end{cases} $

In order to prove that $\frac{1}{\sqrt{n}}X$ and $P$ are asymptotically free we will compute their mixed moments and show that they (asymptotically) coincide with the
corresponding mixed moments of two free random variables $a$ and $b$ with the same distributions as $\frac{1}{\sqrt{n}}X$ and $P$ respectively.

Let $A_1 A_2 \cdots A_m$ be a word in letters $a$ and $b$ (i.e.~for each $i$ either $A_i=a$ or $A_i=b$). %$\frac{1}{\sqrt{n}}X$ and $P$ (i.e.~for each $i$ either $A_i=\frac{1}{\sqrt{n}}X$ or $A_i=P$).
As we are interested only in the trace of $A_1\cdots A_m$ we can without loss of generality assume that the last element $A_m$ of the tuple
is equal to $a$.

As $b$ is a projection we can assume without loss of generality that $b$ does not take neighbouring positions in the tuple $A_1,\dots,A_m$ (i.e.~if $A_i=b$ then $A_{i+1}=a$).

We will need the mixed moments of free $a$ and $b$ mentioned above in order to have something to compare the mixed moments of $\frac{1}{\sqrt{n}}X$ and $P$ to.

It is known that $\varphi(A_1 A_2 \cdots A_m)=\sum_{\pi\in NC(k)}C_{\pi}(a)(\tr(b))^{|max\tau|}$ where $k$ is a number of $a$ in a tuple, 
$\pi$ joins only $a$ and $max\tau$ is a maximal partition of $\{1,2,\dots, m-k\}$ such that $\pi\cup max\tau \in NC(m)$. %sprawdzić to! 

Let $k$ be the number of $i$ such that $A_i=\frac{1}{\sqrt{n}}X$ and let $B_1, B_2, \dots, B_k$ be the same word as $A_1,\dots,A_m$ but in letters $\frac{1}{\sqrt{n}}X$ and $P\frac{1}{\sqrt{n}}X$ (i.e. for every
$l$ either $B_l=\frac{1}{\sqrt{n}}X$ or $B_l=P\frac{1}{\sqrt{n}}X$ and the product $B_1\cdots B_k$ is equal to the product $A_1\cdots A_m$).

It is easy to check that $PX$ is the matrix 
\begin{equation} \label{eq:macierzjucysia}
X=\left[
\begin{matrix}  
 0& 1 & 1 & 1 & \dots & 1 & 1\\
 1& 0 & (1,2) & (1,3) & \dots & (1,q-1)& (1,q)  \\
 1& (1,2) & 0 & (2,3) & \dots  & (2,q-1) & (2,q) \\
 1& (1,3) & (2,3) & 0  & \dots  & (3,q-1) & (3,q)\\
 \vdots & \vdots & \vdots   & \vdots & \ddots & \vdots & \vdots\\
 0& 0 & 0 & 0 & \dots  & 0 & 0 \\
 \end{matrix} \right],
 \end{equation}
  where the last $n-k$ rows consist only of zeros.

We shall now define the \emph{Kreweras complementation map} $K$ from $NC(n)$ to $NC(n)$.
Let $\pi$ be a noncrossing partition of the set $\{1,2,\dots,n\}$. Between points $1,2,\dots,n$ insert new points
$1',2',\dots,n'$ in the following way: $1,1',2,2',\dots,n,n'$. Draw all blocks of $\pi$ and then draw a maximal partition $\pi'$ of $1',2',\dots,n'$
such that $\pi\cup\pi'$ is a non-crossing partition of $1,1',2,2',\dots,n,n'$. Such $\pi'$ is called a Kreweras complement of $\pi$ and will be denoted by
$K(\pi)$. 
%(see the picture).

%\includegraphics{kreweras.jpg}

\begin{lem}[Bia98, Theorem 1.3]
For all $A>1$ and $m$ positive integer, there exists a constant $K>0$ such that, for all $A-balanced$ Young diagrams $\lambda$,
and all permutations $\sigma \in S_{|\lambda|}$ satisfying $|\sigma|\leq m$, one has
$$
|\tr \rho_{\lambda}(\sigma)-\prod_{c|\sigma}|\lambda|^{-|c|-1}C_{|c|+2}(\lambda)|    \leq      K |\lambda |^{-1-\frac{|\sigma|}{2}},
$$
where the product is over the disjoint cycles of the permutation $\sigma$.
\end{lem}

The following Lemma is a reformulation of Theorem 1.3 from [Bia98].

\begin{lem} \label{lemat5}
%The non-commutative probability space $(\mathcal{A},\varphi)$ we work with depend on $n$. Let $n$ go to infinity.
Let balanced Young diagrams $\lambda_1$ and $\lambda_2$ corresponding to $\rho_1$ and $\rho_2$ in the definition of $\varphi$ have, in the limit when
$n$ goes to infinity, some limit shapes
$\Lambda_1$ and $\Lambda_2$.
Let $\sigma$ be a product of some tuple
of Jucys-Murphy transpositions satisfying $\pi \approx (a_1,\dots,a_m) \sim (A_1,\dots,A_m)$ and assume that $\sigma \in S_n \times S_n \times \{e\}$.
Let $\sigma_1, \sigma_2$
be such that $\operatorname{supp}(\sigma_1)\subset\{1,\dots,n\}$, $\operatorname{supp}(\sigma_2)\subset\{n+1,\dots,2n\}$ and $\sigma=\sigma_1 \sigma_2$
where $\operatorname{supp}$ denotes the support of a permutation.

Then $$n^{\frac{|\sigma|}{2}} \varphi(a_1\dots a_m) \rightarrow \prod_{c|\sigma_1}C_{|c|+2}^{\mu_{\Lambda_1}}\prod_{c|\sigma_2}C_{|c|+2}^{\mu_{\Lambda_2}},$$
where $C^{\mu_\Lambda}_k$ denotes the $k$-th free cumulant of $\mu_{\Lambda}$.  
\end{lem}

Asymptotic behaviour of characters of symmetric groups was given in [BiaREF].
\begin{lem}
Let $\lambda_n$ be a sequence of $C$-balanced Young diagrams and $\rho_n$ the corresponding representations of $S_n$. Fix a permutation $\sigma\in S_k$ and note
that $\sigma$ can be treated as an element of $S_n$ if we add $n-k$ additional fixpoints. There exists a constant $K$ such that
\[
\big{|}  \tr(\rho(\sigma))  \big{|} \leq  K n^{\frac{-|\sigma|}{2}}.
\]
\end{lem}

\begin{defi} 
Let $(n)_{k}= n(n-1)\cdots(n-k+1)$ denote the product of descending integers.
\end{defi}

The following computation of the moments of $\frac{1}{\sqrt{n}}X$ was carried out by Biane in [Bia98 proof of Prop.~3.3] with a difference that Biane's
random variable was not normalized by $\frac{1}{\sqrt{n}}$.
\begin{multline*}
\varphi(X^k)=n^{-\frac{k}{2}}\tr\rho(\tr X^k)=\\
\frac{n^{-\frac{k}{2}}}{n+1}
 \sum_{0\leq i_1\neq i_2 \neq \cdots \neq i_n\neq i_1 \leq n} \tr\rho((i_1,i_2)(i_2,i_3)\cdots(i_n,i_1))=(\diamond)
\end{multline*}
Biane has the following way of dealing with the above sum: \emph{`We shall decompose the set of $n$-tuples $(i_1,i_2,\dots,i_n)$ occuring in the above sum according to the set $J$ of places $r_1<r_2<\dots<r_k$
such that $i_{r_j}=0$. For each $J\subset\{1,2,\dots,n\}$ and $i_1,i_2,\dots,i_n$ such that $J=\{l:i_l=0\}$ let $\pi$ be the partition of 
$\{1,2,\dots,n\}\backslash J$ induced by $i_1,i_2,\dots,i_n$, namely $j$ and $k$ belong to the same component of $\pi$ if and only if $i_j=i_k\neq0$.
Clearly the conjugacy class of $(i_1 i_2)(i_2 i_3)\cdots(i_n i_1)$ in $S_q$ depends only on $J$ and $\pi$. We shall denote by $h(\pi)$ this conjugacy class,
and by $|h(\pi)|$ the length of any permutation belonging to it.'} 
$$
(\diamond)=\frac{n^{-\frac{k}{2}}}{n+1} \sum_{J\subset\{1,2,\dots,k\}} \sum_{\pi\in P_a(J,k)} (n)_{|\pi|}
\tr\rho(h(\pi))= (\star)
$$
where $P_a(J,k)$ is the set of all admissible partitions of $\{1,2,\dots,n\}\backslash J$, i.e.~such that $i$ and $i+1$ never
belong to the same component of $\pi$ (we make a convention that if $i=n$ then $i+1=1$). The following Lemmas was proved by Biane:

\begin{lem}[Bia98, Lemma4.3.1]
If $J=\emptyset$ and $\pi$ has a crossing, then $|h(\pi)|\geq 2|\pi|-n$.
\end{lem}

\begin{lem}[Bia98, Lemma4.3.2]
If $J\neq \emptyset$, then $|h(\pi)|\geq 2|\pi|-n$.
\end{lem}

\begin{lem}[Bia98, Lemma4.3.3]
The cycles of any permutation in $h(\pi)$ are in one-to-one correspondence with blocks of $K(\pi)$ and the
order of a cycle is less by one than number of elements of the corresponding block.
\end{lem}

Using Lemmas 3 and 5 we get:

$$
(\star)=\frac{n^{-\frac{k}{2}}}{n+1} \sum_{\pi\in P_a(k)} (n)_{|\pi|}
\tr\rho(h(\pi))= 
$$

Now from Lemmas 3 and 4 we have

$$
=\sum_{\pi\in P_a(k)\cap NC(k)} \frac{(n)_{|\pi|} n^{-\frac{k}{2}}}{n+1} 
\tr\rho(h(\pi)) + o(1)= 
$$
$$
=\sum_{\pi\in P_a(k)\cap NC(k)} \frac{(n)_{|\pi|} n^{-\frac{k}{2}}}{n+1} n^{-\frac{|h(\pi)|}{2}} n^{\frac{|h(\pi)|}{2}}
\tr\rho(h(\pi))+ o(1)= 
$$

From Lemma 2 we get
$$
=\sum_{K(\pi)\in NC_{>1}(k)} \frac{(n)_{|\pi|} n^{-\frac{k}{2}}}{n+1} n^{-\frac{|h(\pi)|}{2}} C_{K(\pi)(\mu_{\lambda})}+ o(1)=(\heartsuit)
$$
\begin{lem}
The length of any permutation in $h(\pi)$ is equal $k-2|K(\pi)|$. 
\end{lem}
$$
(\heartsuit)=\sum_{K(\pi)\in NC_{>1}(k)} \underbrace{\frac{(n)_{|\pi|} n^{-\frac{k}{2}}}{n+1} n^{-\frac{|h(\pi)|}{2}}}_{\text{this tends to 1.}} C_{K(\pi)(\mu_{\lambda})}\rightarrow
$$

$$
\rightarrow\sum_{K(\pi)\in NC_{>1}(k)} C_{K(\pi)(\mu_{\lambda})}.
$$

Let us now compute the mixed moment of $\frac{1}{\sqrt{n}}X$ and $P$ by repeating Biane's computation:

$$
\varphi(A_1 A_2 \cdots A_m)=\varphi(B_1 B_2 \cdots B_k)=
$$
$$
=\frac{n^{-\frac{k}{2}}}{n+1} \sum_{0\leq i_1\neq i_2 \neq \cdots \neq i_n\neq i_1 \leq n } \tr\rho((i_1,i_2)(i_2,i_3)\cdots(i_n,i_1))=
$$
where indexes $i_j$ such that $B_j=P\frac{1}{n}X$ are bounded by $k$.
$$
=\frac{n^{-\frac{k}{2}}}{n+1} \sum_{J\subset\{1,2,\dots,k\}} \sum_{\pi\in P_a(J,k)} \frac{(\Tr P)_S (n-S)_{|\pi|-S}}{(n)_{|\pi|}}(n)_{|\pi|}
\tr\rho(h(\pi)),
$$
where $S$ is the number of blocks $b$ of $\pi$ such that there exists $i\in b$ such that $B_i=P\frac{1}{n}X$.
Now the only difference between the $k$-th moment of $X$ and the above formula is the factor $\frac{(\Tr P)_S (n-S)_{|\pi|-S}}{(n)_{|\pi|}}$.

\begin{lem}
$\frac{(\Tr P)_S (n-S)_{|\pi|-S}}{(n)_{|\pi|}}\rightarrow (\tr P)^{|\max\tau |}$.
\end{lem}
We leave the proof as an excercise for the reader.

By repreating the computation of the $k$-th moment of $X$ we obtain

$$
\sum_{K(\pi)\in NC_{>1}(k)} (\tr P)^{\max\tau} C_{K(\pi)(\mu_{\lambda})}+o(1).
$$
which proves that $\frac{1}{\sqrt{n}}X$ and $P$ are asymptotically free.

We can replace $P$ with $Q$ defined as follows: 
$$Q(\tau) = \begin{cases} \tau &\mbox{if } \tau(q+1)\in\{1,2,\dots,k,q+1\} \\
0 & \mbox{if } \tau(q+1)\in\{k+1,k+2,\dots,q \}.  \end{cases} $$

Such a $Q$ is represented by a matrix 
\begin{equation} \label{eq:macierzjucysia}
P=\left[
\begin{matrix}   1 & 0 & 0 & \dots & 0 & 0 \\
 0 & Q_1 & 0 & \dots  & 0 & 0 \\
0 & 0 & Q_2 & \dots  & 0 & 0 \\
 \vdots & \vdots & \vdots & \ddots  & \vdots & \vdots \\
 0 & 0 & 0 & \dots  & Q_{q-1} & 0 \\
 0 & 0 & 0 & \dots & 0 & Q_q
 \end{matrix} \right],
 \end{equation}
where $Q_j(\sigma)=\begin{cases} \sigma &\mbox{if } \sigma(j)\in\{1,2,\dots,k,q+1\} \\
0 & \mbox{if } \sigma(j)\in\{k+1,k+2,\dots,q \}.  \end{cases} $

%The matrix of $QX$ is then 
%\begin{equation} \label{eq:macierzjucysia}
%X=\left[
%\begin{matrix}  
% 0& 1 & 1 & 1 & \dots & 1 & 1\\
% 1& 0 & (1,2) & (1,3) & \dots & 0&0 \\
% 1& (1,2) & 0 & (2,3) & \dots  & 0 & 0 \\
% 1& (1,3) & (2,3) & 0  & \dots  & 0 & 0 \\
 %\vdots & \vdots & \vdots   & \vdots & \ddots & \vdots & \vdots\\
% 1& (1,q) & (2,q) & (3,q) & \dots  & 0 & 0 \\
% \end{matrix} \right],
 %\end{equation}
%where the last $n-k$ columns consists of zeroes except for the first row which consist only ones.
%As the matrix $X$ is symmetric we can prove that $X$ and $Q$ are free in the similar way we did it for $X$ and $P$.
But if we change the identification map $f: S_q\times\{1,2,\dots,q\}\rightarrow S_{q+1}$ from $f(\sigma,\tau)\mapsto\sigma\tau$ to
$f(\sigma,\tau)\mapsto\tau\sigma$ then the matrix of the right multiplication by $X$ is the same as the matrix of
the left multiplication in the previous identification with a difference that the entries are acting by the left regular representation.
It is easy to check that in this new language the matrix of $Q$ is equal to the matrix of $P$ in the old language and the same
proof gives us the freeness of $X$ and $Q$.

\begin{center}

\textbf{REFERENCES}

\end{center}

\medskip

\noindent [Bia98] Philippe Biane. Representations of symmetric groups and free probability. 
 								  Adv. Math., 138(1):126–181, 1998.
\medskip

\noindent [Bia01a] Philippe Biane. Approximate factorization and concentration for characters of
                   symmetric groups. Internat. Math. Res. Notices, (4):179–192, 2001. 

\medskip

\end{document}